\documentclass[12pt,reqno]{amsart}
\usepackage{amsmath,amssymb,amsfonts,amscd,latexsym,amsthm,mathrsfs}
\usepackage[usenames]{color}
\usepackage[unicode]{hyperref}
\renewcommand{\d}{{\mathrm d}}

\theoremstyle{remark}

\begin{document}

\title{The birthday boy problem}

\date{14 August 2021. \emph{Updated}: 1 April 2023}

\author{Wadim Zudilin}
\address{Department of Mathematics, IMAPP, Radboud University, PO Box 9010, 6500~GL Nijmegen, Netherlands}
\email{w.zudilin@math.ru.nl}


\begin{abstract}
In their recent (2021) preprint, Robert Dougherty-Bliss, Christoph Koutschan and Doron Zeilberger come up with a powerful strategy to prove the irrationality, in a quantitative form, of some numbers that are given as multiple integrals or quotients of such.
What is really missing there, for many examples given, is an explicit identification of those irrational numbers;
the authors comment on this task, ``The output file [\dots] contains many such conjectured evaluations, (very possibly many of them are equivalent via a
hypergeometric transformation rule) and we challenge [\dots], the birthday boy, or anyone
else, to prove them.''
Without an identification, the numbers are hardly appealing to human (number theorists).
The goal of this note is to outline a strategy to do the job and illustrate it on several promising entries discussed in the preprint cited above.
\end{abstract}

\maketitle

Recently, R.~Dougherty-Bliss, Ch.~Koutschan and D.~Zeilberger \cite{D-BKZ21} came up with a large collection of numbers, for which the irrationality can be established using explicit constructions of rational approximations.
The latter are based on generalisations of the Beukers famous integrals \cite{Be79} for $\zeta(2)$ and $\zeta(3)$,
namely, on the following double and triple integrals:
\begin{align*}
I_2(a_1,a_2,b_1,b_2;n)
&=\frac{\Gamma(2-a_1-a_2)\,\Gamma(2-a_1-a_2)}{\Gamma(1-a_1)\,\Gamma(1-a_2)\,\Gamma(1-b_1)\,\Gamma(1-b_2)}
\\ &\quad\times
\iint\limits_{(0,1)^2}\frac{x^{n-a_1}(1-x)^{n-a_2}y^{n-b_1}(1-y)^{n-b_2}}
{(1-xy)^{n+1}}\,\d x\,d y
\end{align*}
and
\begin{align*}
J_3(a,b,c,d,e;n)
&=\iiint\limits_{(0,1)^3}\frac{x^{n+b}(1-x)^{n+c}y^{n+e}(1-y)^{n+a}z^{n+a}(1-z)^{n+c}}
{(1-(1-xy)z)^{n+d}}\,\d x\,d y\,\d z.
\end{align*}
The main results in \cite{D-BKZ21} concern the irrationality of the quantities
$$
C(a_1,a_2,b_1,b_2)=I_2(a_1,a_2,b_1,b_2;0)
\quad\text{and}\quad
K(a,b,c,d,e)=\frac{J_3(a,b,c,d+1,e;0)}{J_3(a,b,c,d,e;0)}
$$
together with their irrationality measures, at the moment experimentally but with plausible expectations that those observations can be converted into rigour.
A side consequence of their findings, completely rigorous(!), is presence of the Ap\'ery-like second-order recurrence equations, hence nice and regular forms of irregular continued fractions for the numbers in question (also for those whose irrationality is not yet established).

But what are these numbers? For $C(a_1,a_2,b_1,b_2)$, an explicit identification is already provided \cite{D-BKZ21}:
$$
C(a_1,a_2,b_1,b_2)
={}_3F_2\bigg(\begin{matrix} 1, \, 1-a_1, \, 1-b_1 \\ 2-a_1-a_2, \, 2-b_1-b_2 \end{matrix} \biggm| 1\bigg).
$$
Here and in what follows we use the standard notation for $_{p+1}F_p$-hypergeometric functions \cite{Ba35}.
In certain cases, the $_3F_2$-hypergeometric expression can be evaluated in terms of gamma functions only;
for example \cite{D-BKZ21},
$$
C(0,-1/2,1/6,-1/2)=-24-\frac{\sqrt\pi\,\Gamma(7/3)}{\Gamma(-1/6)}
$$
and
$$
C(-2/3,-1/2,1/2,-1/2)=\frac{13}2-\frac{6\Gamma(19/6)}{\sqrt\pi\,\Gamma(8/3)},
$$
and the second (irrational) summands are rationally proportional to $C_1=2^{-1/3}\sqrt3\,\*\Gamma(1/3)^3/\pi$
and $C_2=2^{-1/3}\Gamma(1/3)^3/\pi^2$, respectively.
Both are known to be transcendental (because of the algebraic independence of $\Gamma(1/3)$ and $\pi$), however the estimates for their irrationality measures seem to be new.
The constructions in \cite{D-BKZ21} also result in irregular continued fractions for $C_1$ and $C_2$.

What can be said about $K(a,b,c,d,e)$? Apart from a few cases, an identification of these numbers with some `recognisable' expressions is not developed in \cite{D-BKZ21}, mainly, because of miss of a human touch in executing the computer algebra systems, \texttt{Maple} and \texttt{Wolfram Mathematica}.
Citing \cite{D-BKZ21},
``The output file [\dots] contains many such conjectured evaluations, (very possibly many of them are equivalent via a hypergeometric transformation rule) and we challenge [\dots], the birthday boy, or anyone else, to prove them.''
The missing identification is crucial in appreciating the power of the results in \cite{D-BKZ21}!
For this reason, we indicate below a recipe to perform the task and illustrate it on several examples.

The principal ingredient is a theorem of Nesterenko (see \cite[Proposition 1]{Zu04} for the form adopted to present needs):
\begin{align*}
J_3&=J_3(a,b,c,d,e;0)
=\iiint\limits_{[0,1]^3}\frac{x^b(1-x)^cy^e(1-y)^az^a(1-z)^c}{(1-z+xyz)^d}\,\d x\,\d y\,\d z
\\
&=\frac{\Gamma(1+a)\Gamma(1+c)}{\Gamma(d)\Gamma(2+a+c-d)}
\\ &\kern-5mm\times
\frac1{2\pi i}\int\limits_{t_1-i\infty}^{t_1+i\infty}
\frac{\small\Gamma(d+t)\Gamma(1+a+t)\Gamma(1+b+t)\Gamma(1+e+t)\cdot\Gamma(-t)\Gamma(1+c-d-t)}{\small\Gamma(2+a+e+t)\Gamma(2+b+c+t)}\,\d t,
\end{align*}
where
$$
a,b,c,e>-1 \quad\text{and}\quad
-\min\{a+1,b+1,d,e+1\}<t_1<c-d+1.
$$
This leads to the Meijer $G$-function expression
$$
J_3(a,b,c,d,e)=\frac{\Gamma(1+a)\Gamma(1+c)}{\Gamma(d)\Gamma(2+a+c-d)}\,
G_{4,4}^{2,4}\bigg(1\biggm|\begin{matrix} 1-d, \, -a, \, -b, \, -e \\
1-d+c, \, 0; \; -1-b-c, \, -1-e-a \end{matrix}\bigg),
$$
and also, in the case $c-d\notin\mathbb Z$, to the hypergeometric reduction
\begin{align*}
J_3&=\frac{\pi}{\sin\pi(d-c)}
\\ &\;\times
\biggl(
\frac{\Gamma(a+1)^2\Gamma(b+1)\Gamma(c+1)\Gamma(e+1)}
{\Gamma(d-c)\Gamma(a+e+2)\Gamma(b+c+2)\Gamma(a+c-d+2)}
\\ &\;\qquad\quad\times
{}_4F_3\biggl(\begin{matrix} a+1, \, b+1, \, d, \, e+1 \\ d-c, \, a+e+2, \, b+c+2 \end{matrix} \biggm|1\biggr)
\\ &\;\qquad
-\frac{\Gamma(a+1)\Gamma(c+1)\Gamma(b+c-d+2)\Gamma(c+1)\Gamma(e+c-d+2)}
{\Gamma(c-d+2)\Gamma(d)\Gamma(a+e+c-d+3)\Gamma(b+2c-d+3)}
\\ &\;\qquad\qquad\quad\times
{}_4F_3\biggl(\begin{matrix} a+c-d+2, \, b+c-d+2, \, c+1, \, e+c-d+2 \\ c-d+2, \, a+e+c-d+3, \, b+2c-d+3 \end{matrix} \biggm|1\biggr) \biggr).
\end{align*}

Equipped with this formula and computer algebra systems we have the following findings:
\begin{alignat*}{2}
K(0,0,0,2/3,1/3)
& 
=-\frac{K_1-2}{2(K_1-3)},
&\quad\text{where}\;
K_1&=\log3+\frac{\pi}{\sqrt3},
\displaybreak[2]\\
K(0,0,0,1/3,2/3)
& 
=-\frac{2(K_2+1)}{K_2+1/2},
&\quad\text{where}\;
K_2&=\log3+\frac{\pi}{\sqrt3},
\displaybreak[2]\\
K(0, 1/3, 2/3, 1/3, 2/3)
& 
=-\frac{20(7-54K_3)}{52-405K_3},
&\quad\text{where}\;
K_3&=\frac{\Gamma(2/3)^3}{\Gamma(1/3)^3},
\displaybreak[2]\\
K(0, 1/5, 0, 3/5, 2/5)
& 
=-\frac{4(1-4K_4)}{5-24K_4},
&\quad\text{where}\;
K_4&=\frac1{\sqrt5}\,\log\frac{\sqrt5+1}2,
\displaybreak[2]\\
K(1/7, 0, 2/7, 3/7, 4/7)
& 
=-\frac{189(8-5K_5)}{832 - 525K_5},
&\quad\text{where}\;
K_5&=\frac{2^{2/7}\sqrt\pi\,\Gamma(9/14)}{\cos(3\pi/14)\,\Gamma(4/7)^2},
\end{alignat*}
and similarly looking evaluations for 
$K(1/7, 0, 2/7, 5/7, 3/7)$, $K(1/7, 0, 3/7, 4/7, 5/7)$ and $K(1/7, 0, 4/7, 2/7, 5/7)$,
while the expressions found for $K(2/7, 0, 3/7, 4/7, 5/7)$, $K(0, 1/2, 1/2, 1/3, 1/6)$ and $K(0, 1/2, 1/2, 1/6, 1/3)$ involve $_3F_2$-hypergeometric functions that seem to be not reducible to a simpler form.
Finally,
$$
K(1/3, 0, 2/3, 1/2, 5/6)
=-\frac{25 (84 - 53\cdot 2^{2/3})}{9(216-137\cdot 2^{2/3})}.
$$
Perhaps, a real pearl in this collection of `quantitatively' irrational numbers (admitting regular-form irregular continued fractions) is the number~$K_3$.

\medskip\noindent
\textbf{Acknowledgements}.
I am thankful to Christoph Koutschan for related discussions and corrections.

\end{document}